# MINIMAL $f^Q$-MARTINGALE MEASURES FOR EXPONENTIAL LÉVY PROCESSES

By Monique Jeanblanc, Susanne Klöppel[1] and Yoshio Miyahara[2]

*Université d'Evry and Institut Europlace de France, Vienna University of Technology and Nagoya City University*

Let $L$ be a multidimensional Lévy process under $P$ in its own filtration. The $f^q$-minimal martingale measure $Q_q$ is defined as that equivalent local martingale measure for $\mathcal{E}(L)$ which minimizes the $f^q$-divergence $E[(dQ/dP)^q]$ for fixed $q \in (-\infty, 0) \cup (1, \infty)$. We give necessary and sufficient conditions for the existence of $Q_q$ and an explicit formula for its density. For $q = 2$, we relate the sufficient conditions to the structure condition and discuss when the former are also necessary. Moreover, we show that $Q_q$ converges for $q \searrow 1$ in entropy to the minimal entropy martingale measure.

**1. Introduction.** Lévy models are very popular in finance due to their tractability and their good fitting properties. However, Lévy models typically yield incomplete markets. This raises the question of which measure one should choose for valuation or pricing of nonhedgeable payoffs. Very often, a measure is chosen which minimizes a particular functional over the set $\mathcal{M}^e(S)$ of equivalent local martingale measures for the underlying assets $S$. This choice can be motivated by a dual formulation of a primal utility maximization problem; see Kramkov and Schachermayer [15] and Frittelli [7]. If $P$ denotes the subjective measure, then the functional on $\mathcal{M}^e(S)$ is typically of the form $Q \mapsto E_P[f(dQ/dP)]$, where $f$ is a convex function on $(0, \infty)$. Then $f(Q|P) := E_P[f(dQ/dP)]$, known as the *f-divergence* of $Q$ with respect to $P$, is a measure for the distance between $Q$ and $P$; see

Received November 2006; revised March 2007.
[1]Supported by the NCCR FINRISK and the Christian Doppler Research Association.
[2]Supported by Grant-in-Aid for Scientific Research No.16540113, JSPS.

Liese and Vajda [18] for a textbook account. Hence, one chooses as pricing measure the martingale measure which is closest to $P$ with respect to some $f$-divergence.

In this article, we consider $f^q(Q|P)$ corresponding to $f^q(z) = z^q$ for $q \in (-\infty, 0) \cup (1, \infty)$. The optimal measure $Q_q$ is then called the $f^q$-*minimal martingale measure*. More precisely, we work on a probability space $(\Omega, \mathcal{F}, P)$ equipped with a filtration which is the $P$-augmentation of that generated by a $d$-dimensional Lévy process $L$ and model the traded assets $S$ as the stochastic exponential $S = \mathcal{E}(L)$. Based on an explicit formula for $f^q(Q|P)$ in terms of the Girsanov parameters $(\beta, Y)$ of $Q$, we show that $f^q(Q|P)$ is reduced if $Q$ is replaced by some $\overline{Q} \in \mathcal{M}^e(S)$, which is defined via its Girsanov parameters $(\overline{\beta}, \overline{Y})$; see Theorem 2.6. The measure $\overline{Q}$ preserves the Lévy property of $L$, that is, $L$ is a $\overline{Q}$-Lévy process. We deduce that $Q_q$ also has this property and that for minimization of $f^q(Q|P)$, it suffices to consider those $Q \in \mathcal{M}^e(S)$ which preserve the Lévy property of $L$. As illustrated in Theorem 2.7, this allows the reduction of the minimization of $f^q(Q|P)$ to a deterministic convex optimization problem $(\mathcal{P}_q)$ whose solution corresponds to the Girsanov parameters of $Q_q$ and hence provides an explicit formula for $Q_q$. In particular, $(\mathcal{P}_q)$ has a solution if and only if $Q_q$ exists. By formally applying the Kuhn–Tucker theorem to $(\mathcal{P}_q)$, we obtain, in Theorem 2.9, conditions $(\mathcal{C}_q)$, which are sufficient for the existence of a solution to $(\mathcal{P}_q)$ and which can easily be verified in practice. From $(\mathcal{C}_q)$, one can immediately deduce the solution to $(\mathcal{P}_q)$ and hence obtain the explicit formula for $Q_q$. For $q = 2$, we relate $Q_2$ to the *variance optimal signed martingale measure* and $(\mathcal{C}_2)$ to the *structure condition* (SC). In Theorem 3.1, we specify some cases in which $(\mathcal{C}_2)$ is not only sufficient, but also necessary for the existence of $Q_2$; this requires an additional integrability condition for $L$. Finally, we prove that, under some technical assumptions, $Q_q$ converges for $q \searrow 1$ in entropy to the *minimal entropy martingale measure* $P_e$ and also that the corresponding Girsanov parameters converge; $P_e$ is defined as that measure $Q$ which minimizes the divergence corresponding to $f(z) := z \log z$ over all local martingale measures for $S$. The convergence is shown by an application of the implicit function theorem; for this, $(\mathcal{C}_q)$ must be rewritten as the root of an appropriate function. A concluding example illustrates that the technical conditions we must impose for convergence can all be satisfied in a reasonable model.

Some of the results and concepts have been studied before. In [6], Esche and Schweizer use an approach similar to our Theorem 2.6, but for $P_e$ instead of $Q_q$. However, we define $(\overline{\beta}, \overline{Y})$ above by a pointwise criterion, whereas they obtain them from averaging the Girsanov parameters of the original measure $Q$ with respect to an appropriate measure on $\Omega \times [0, T]$. Therefore, they require an extensive approximation procedure in order to prove that $\overline{Q}$ is, in fact, a martingale measure. We thus do not only generalize their approach to $Q_q$, but, more importantly, significantly simplify it. Although problem $(\mathcal{P}_q)$ itself seems not to have been studied before, the sufficient conditions $(\mathcal{C}_q)$ for $q \in (-\infty, 0)$ also appear in Kallsen [14] and Goll and Rüschendorf [10]. However, they do only state $(\mathcal{C}_q)$ without motivating its definition, which is done in our approach. Very recently, we became familiar with an article by Choulli, Stricker and Li [4] in which they obtain $(C_q)$ for general $q$, but for the *minimal Hellinger martingale measure* $\mathcal{Q}$ of order $q$ instead of $Q_q$. In Section 2, we show that $Q_q$ and $\mathcal{Q}$ coincide in our Lévy setting. Convergence of $Q_q$ to $P_e$ for $q \searrow 1$ was studied for continuous processes by Grandits and Rheinländer in [11], but an extension to processes which also have a jump part seems to be missing. The discussion of the conditions when $(\mathcal{C}_2)$ is even necessary for the existence of $Q_2$ also seems to be new. Finally, we think it is remarkable that our approach is very intuitive and accessible to anyone interested in risk-neutral measures.

The paper is structured as follows. In Section 2, we give necessary and sufficient conditions for the existence of $Q_q$ and an explicit formula for its density. Section 3 covers additional results for the special case $q = 2$. Convergence of $Q_q$ to $P_e$ is presented in Section 4. Appendix A contains the required results on Lévy processes, Appendix B some auxiliary calculations and Appendix C some proofs omitted from the main body the article.

**2. Structure and existence of the $f^q$-minimal martingale measure.** In this section, we give necessary and sufficient conditions for the existence of the $f^q$-minimal martingale measure (qMMM) $Q_q$ for an exponential Lévy process $S$; in the entire section, $q \in I := (-\infty, 0) \cup (1, \infty)$ is arbitrary but fixed. In particular, we provide explicit formulas for the density of $Q_q$.

To begin with, some notation and conventions are introduced; our basic reference is Jacod and Shiryaev [13]. Throughout the article, we work on a filtered probability space $(\Omega, \mathcal{F}, \mathbb{F}, P)$, where $\mathbb{F} = \mathbb{F}^L$ is the $P$-augmentation of the filtration generated by a $d$-dimensional Lévy process $L = (L_t)_{0 \le t \le T}$ with characteristic triplet $(b, c, K)$ with respect to the truncation function $h(x) := x\mathbf{I}_{\{\|x\| \le 1\}}$ and where $T$ is a finite time horizon. If it exists, we choose a right-continuous version of any process. The random measure associated with the jumps of $L$ is denoted by $\mu^L$ and $\nu^P(dx, dt) = K(dx)\, dt$ is the predictable $P$-compensator of $\mu^L$; all required background on Lévy processes and unexplained terminology can be found in Appendix A. The stochastic exponential of $L$ is denoted by $S := \mathcal{E}(L) = (\mathcal{E}(L^1), \ldots, \mathcal{E}(L^d))^*$, where $*$ indicates the transpose of a vector. We assume that the process $S$ is strictly positive, that is, that $\Delta L^i > -1$, $P$-a.s. for $i \in \{1, \ldots, d\}$. By $\mathcal{M}^e(S)$, we denote the set of all equivalent local martingale measures for $S$; note that $\mathcal{M}^e(S) = \mathcal{M}^e(L)$ due to Ansel and Stricker [1], Corollary 3.5. For $f^q(z) := z^q$ and $q \in I$, the qMMM $Q_q \in \mathcal{M}^e(S)$ is characterized by the property that it minimizes the $f^q$-divergence

$$f^q(Q|P) := E[f^q(Z_T^Q)] = E[(Z_T^Q)^q]$$

over all $Q \in \mathcal{M}^e(S)$; by $E[\cdot]$, we denote the expectation with respect to $P$ and for any $Q \ll P$, its real-valued density process $Z^Q = (Z_t^Q)_{0 \leq t \leq T}$ with $Z_t^Q := E[dQ/dP|\mathcal{F}_t]$ is defined with respect to $P$. In particular, we require that $f^q(Q_q|P) < \infty$, that is, that $Q_q$ is contained in the set

$$\mathcal{Q}^q := \{Q \in \mathcal{M}^e(S) | f^q(Q|P) < \infty\}.$$

In the present Lévy setting, any $Q \approx P$ can be fully described by its Girsanov parameters $(\beta, Y)$ with respect to $L$ and we write $Q = Q^{(\beta,Y)}$ to emphasize this; see Proposition A.3 in Appendix A. It is well known that the $P$-Lévy process $L$ is also a $Q^{(\beta,Y)}$-Lévy process for some $Q^{(\beta,Y)}$ equivalent to $P$, if and only if $\beta$ and $Y$ are time-independent and deterministic; see Corollary A.5 in Appendix A. The set of all $Q \in \mathcal{M}^e(S)$ having this property will be of importance later and we denote it by

$$\overline{\mathcal{Q}} := \{Q^{(\beta,Y)} \in \mathcal{M}^e(S) | (\beta, Y) \text{ time-independent and deterministic}\}.$$

REMARK 2.1. All results are stated for $S = \mathcal{E}(L)$. However, we could equivalently work with $S = e^{\widetilde{L}}$, where $\widetilde{L}$ is a $P$-Lévy process with $P$-characteristic triplet $(\widetilde{b}, \widetilde{c}, \widetilde{K})$ since $e^{\widetilde{L}} = \mathcal{E}(L)$ if $L$ has characteristic triplet

$$b = \widetilde{b} + \tfrac{1}{2}(\widetilde{c}_{11}, \ldots, \widetilde{c}_{dd})^* + \int_{\mathbb{R}^d} (h(e^x - \mathbf{1}) - h(x))\widetilde{K}(dx),$$

$$c = \widetilde{c},$$

$$K(G) = \int_{\mathbb{R}^d} \mathbf{I}_{\{(e^x - \mathbf{1}) \in G\}} \widetilde{K}(dx),$$

where $\mathbf{1} = (1, \ldots, 1)^*$; this holds by Itô's formula and is stated explicitly in Goll and Kallsen [9], Lemma A.8. Analogously, discounting with respect to some numeraire can also be captured by a modification of the characteristic triplet of $L$; see Corollary A.7 in Goll and Kallsen [9].

2.1. *Reducing the problem.* In this subsection, we simplify the characterization of the qMMM $Q_q \in \mathcal{Q}^q$

(2.1) $$f^q(Q_q|P) = \inf_{Q \in \mathcal{M}^e(S)} f^q(Q|P) < \infty.$$

More precisely, we show that $Q_q$ can be fully described as the solution to a deterministic optimization problem in $\mathbb{R}^d$ which has a solution if and only if $Q_q$ exists. For this, we exploit the fact that the dynamics of the Lévy process $L$, and hence those of $S$, are time-independent; this becomes apparent in $(b, c, K)$. In fact, we show that it does not only suffice to minimize the $f^q$-divergence $f^q(Q|P)$ over the set $\mathcal{Q}^q$, but over those $Q = Q^{(\beta, Y)}$ with deterministic, time-independent Girsanov parameters, that is, those which are contained in $\overline{\mathcal{Q}}$.

As a first step in this direction, we give a nice formula for $f^q(Q|P)$. For this, we introduce the strictly convex, nonnegative function

$$g_q(y) := y^q - 1 - q(y-1), \qquad y \in (0, \infty).$$

PROPOSITION 2.2. *Let $Q = Q^{(\beta, Y)} \in \mathcal{Q}^q$ with density process $Z = Z^Q = \mathcal{E}(N)$. The canonical $P$-decomposition of the $P$-submartingale $f^q(Z) = Z^q = 1 + M + A$ is*

$$M = \int Z_-^q \, d\widehat{M} = \int f^q(Z_-) \, d\widehat{M} \qquad \text{with}$$

$$\widehat{M} := qN + g_q(Y) * (\mu^L - \nu^P) = qN^c + (Y^q - 1) * (\mu^L - \nu^P)$$

*and*

$$A = \int Z_-^q \, d\widehat{A} = \int f^q(Z_-) \, d\widehat{A} \qquad \text{with}$$

$$\widehat{A} := \frac{q(q-1)}{2} \langle N^c \rangle + g_q(Y) * \nu^P.$$

*Its multiplicative decomposition is $f^q(Z) = \mathcal{E}(\widehat{M})\mathcal{E}(\widehat{A})$, where $\mathcal{E}(\widehat{M})$ is a strictly positive uniformly integrable $P$-martingale and $\mathcal{E}(\widehat{M})_0 = 1$. With $\frac{dR_q}{dP} := \mathcal{E}(\widehat{M})_T$, the process $\mathcal{E}(\widehat{A}) = e^{\widehat{A}}$ is increasing and $R_q$-integrable and*

$$(2.2) \qquad f^q(Q|P) = E_{R_q}[\mathcal{E}(\widehat{A})_T] = E_{R_q}\left[\exp\left(\int_0^T k_q(\beta_t, Y_t) \, dt\right)\right],$$

*where $k_q(\beta_t, Y_t) := \frac{q(q-1)}{2} \beta_t^* c \beta_t + \int_{\mathbb{R}^d} g_q(Y_t(x)) K(dx)$.*

PROOF. See Appendix C. □

The distinctive feature of Proposition 2.2 is formula (2.2). Note that for $Q \in \overline{\mathcal{Q}} \cap \mathcal{Q}^q$, the dependence on $R_q$ vanishes, and that the formula reduces to the following expression.

COROLLARY 2.3. *If, in Proposition 2.2, we have $Q^{(\beta, Y)} \in \mathcal{Q}^q \cap \overline{\mathcal{Q}}$, then*

$$f^q(Q|P) = e^{T k_q(\beta, Y)} = \exp\left(T\left(\frac{q(q-1)}{2} \beta^* c \beta + \int_{\mathbb{R}^d} g_q(Y(x)) K(dx)\right)\right).$$

PROOF.  This is obvious from the definition of $\overline{\mathcal{Q}}$ and $k$. □

As mentioned above, the Lévy structure is essentially time-independent. This suggests that $Q_q$ should also be of this form, that is, that $Q_q \in \overline{\mathcal{Q}} \cap \mathcal{Q}^q$. We will prove this as follows. By (2.2), we can find, for any $Q^{(\beta,Y)} \in \mathcal{Q}^q$, a pair $(\overline{\omega}, \overline{t}) \in \Omega \times [0,T]$ such that $\overline{\beta} := \beta_{\overline{t}}(\overline{\omega})$ and $\overline{Y}(x) := Y_{\overline{t}}(x, \overline{\omega})$ satisfy

$$e^{Tk_q(\overline{\beta},\overline{Y})} \leq f^q(Q^{(\beta,Y)}|P).$$

However, we must ensure that $(\overline{\beta}, \overline{Y})$ can be identified with the Girsanov parameters of some $\overline{Q} = Q^{(\overline{\beta},\overline{Y})} \in \overline{\mathcal{Q}} \cap \mathcal{Q}^q$ so that, by Corollary 2.3,

$$e^{Tk_q(\overline{\beta},\overline{Y})} = f^q(\overline{Q}|P).$$

This is the content of our next result.

PROPOSITION 2.4.  *Let $\overline{\beta} \in \mathbb{R}^d$ and $\overline{Y}: \mathbb{R}^d \to \mathbb{R}_+$ be a measurable function such that $\overline{Y} > 0$, $K$-a.e.,*

(2.3) $$k_q(\overline{\beta}, \overline{Y}) < \infty$$

*and the* martingale condition *holds, that is,*

$$\int_{\mathbb{R}^d} |x\overline{Y}(x) - h(x)| K(dx) < \infty,$$
(ℳ)
$$b + c\overline{\beta} + \int_{\mathbb{R}^d} (x\overline{Y}(x) - h(x)) K(dx) = 0.$$

*$(\overline{\beta}, \overline{Y})$ are then the Girsanov parameters of some $\overline{Q} := Q^{(\overline{\beta},\overline{Y})} \in \overline{\mathcal{Q}} \cap \mathcal{Q}^q$.*

PROOF.  See Appendix C. □

REMARK 2.5.  Proposition A.7 in Appendix A justifies why we refer to (ℳ) as a martingale condition.

THEOREM 2.6.  *Let $\mathcal{Q}^q \neq \varnothing$. Then:*

1. *for any $Q \in \mathcal{Q}^q$, there exists $\overline{Q} \in \overline{\mathcal{Q}}$ such that*
$$f^q(\overline{Q}|P) \leq f^q(Q|P);$$

2. *for every $\varepsilon > 0$, there exists $\overline{Q} \in \overline{\mathcal{Q}}$ such that*
$$f^q(\overline{Q}|P) \leq \inf_{Q \in \mathcal{M}^e(L)} f^q(Q|P) + \varepsilon;$$

3. *if the qMMM $Q_q$ exists, then $Q_q \in \overline{\mathcal{Q}}$.*

PROOF. 1. By part 1 of Remark A.2, Propositions 2.2 and A.7, and since $Q = Q^{(\beta,Y)} \in \mathcal{Q}^q$, there exists $(\overline{\omega}, \overline{t}) \in \Omega \times [0,T]$ and such that $\overline{\beta} := \beta_{\overline{t}}(\overline{\omega})$ and $\overline{Y}(x) := Y_{\overline{t}}(x, \overline{\omega})$ satisfy $\overline{Y}(x) > 0$, $K$-a.e., ($\mathcal{M}$) and

$$e^{Tk_q(\overline{\beta},\overline{Y})} \leq E_{R_q}\left[\exp\left(\int_0^T k_q(\beta_t, Y_t)\, dt\right)\right] = f^q(Q|P) < \infty;$$

note that $\beta$ and $Y$ satisfy $Y_t(x, \omega) > 0$ $K$-a.e. and ($\mathcal{M}$) for $dP \otimes dt$-a.e. $(\omega, t) \in \Omega \times [0,T]$. The claim then follows from Proposition 2.4 and Corollary 2.3.

2. Since $\mathcal{Q}^q \neq \varnothing$ and by the definition of the infimum, there exists $Q' \in \mathcal{Q}^q$ such that $f^q(Q'|P) \leq \inf_{Q \in \mathcal{M}^e(L)} f^q(Q|P) + \varepsilon$. Thus 2 follows from 1 applied to $Q'$.

3. Since $z \mapsto z^q$ is strictly convex on $\mathbb{R}_+$, $Q_q$ is unique. Thus 3 follows immediately from 1 applied to $Q_q$. $\square$

The proof of Theorem 2.6 is similar to the ansatz used by Esche and Schweizer in [6] for the minimal entropy martingale measure $P_e$; the latter is defined as that local martingale measure $Q$ for $S$ which minimizes the relative entropy $E_Q[\log(dQ/dP)]$. However, for $Q^{(\beta,Y)} \in \mathcal{Q}^q$ we define the Girsanov parameters of a measure $\overline{Q}$ with $f^q(\overline{Q}|P) \leq f^q(Q^{(\beta,Y)}|P)$ by *fixing* $(\overline{\omega}, \overline{t}) \in \Omega \times [0,T]$, that is, we set $\overline{\beta} := \beta_{\overline{t}}(\overline{\omega})$ and $\overline{Y}(x) := Y_{\overline{t}}(x, \overline{\omega})$, whereas they *average* $\beta$ and $Y$ with respect to an appropriate measure on $\Omega \times [0,T]$; in our setting, this would be $dR_q \otimes dt$ with $R_q$ from Proposition 2.2. The advantage of our pointwise approach is that it ensures that $\overline{Q} \in \mathcal{M}^e(S)$, whereas Esche and Schweizer apply Fubini's theorem to prove that $\overline{Q}$ is again a local martingale measure. The latter requires an additional integrability condition on $L$ which is not satisfied in general. Thus they must show that there is a dense subset of local martingale measures with this integrability condition and apply additional approximation procedures. This can all be avoided by our approach.

Next, we state the aforementioned deterministic optimization problem. It uses the fact that Proposition 2.4 provides us with a complete characterization of the set $\overline{\mathcal{Q}} \cap \mathcal{Q}^q$; necessity of (2.3) and ($\mathcal{M}$) follow from Corollary 2.3 and Proposition A.7 in Appendix A. Since we know from Theorem 2.6 that $Q_q$, if it exists, is contained in $\overline{\mathcal{Q}} \cap \mathcal{Q}^q$, we can thus describe $Q_q$ as the solution to the following optimization problem which has a solution if and only if $Q_q$ exists.

($\mathcal{P}_q$): *Find a solution $(\widehat{\beta}_q, \widehat{Y}_q)$ to*

$$minimize \quad k_q(\beta, Y) = \frac{q(q-1)}{2}\beta^* c\beta + \int_{\mathbb{R}^d} g_q(Y(x))K(dx)$$

*over*

$$\mathcal{A}_q := \{(\beta, Y) | \beta \in \mathbb{R}^d, Y : \mathbb{R}^d \to \mathbb{R}_+ \text{ measurable},$$
$$Y > 0, \ K\text{-a.e.}, (\beta, Y) \text{ satisfies } (\mathcal{M}), k_q(\beta, Y) < \infty\}.$$

THEOREM 2.7. *(i) If the qMMM $Q_q = Q^{(\beta_q, Y_q)}$ exists in $\mathcal{Q}^q$, then $(\beta_q, Y_q)$ solves $(\mathcal{P}_q)$, that is, $(\beta_q, Y_q) = (\widehat{\beta}_q, \widehat{Y}_q)$.*

*(ii) If $(\widehat{\beta}_q, \widehat{Y}_q)$ solves $(\mathcal{P}_q)$, then $Q_q$ exists and has Girsanov parameters $(\widehat{\beta}_q, \widehat{Y}_q)$, that is, $Q_q = Q^{(\widehat{\beta}_q, \widehat{Y}_q)}$.*

PROOF. (i) By 3 of Theorem 2.6, we have $Q_q \in \overline{\mathcal{Q}}$ and Proposition A.7 implies that $(\beta_q, Y_q)$ satisfies $(\mathcal{M})$, so $(\beta_q, Y_q) \in \mathcal{A}_q$ by Corollary 2.3. Suppose that there exists $(\beta, Y) \in \mathcal{A}_q$ with $k_q(\beta, Y) < k_q(\beta_q, Y_q) < \infty$. Then, by Proposition 2.4 and Corollary 2.3, there exists $Q := Q^{(\beta, Y)} \in \overline{\mathcal{Q}} \cap \mathcal{Q}^q$ such that

$$f^q(Q|P) = e^{Tk_q(\beta, Y)} < e^{Tk_q(\beta_q, Y_q)} = f^q(Q_q|P),$$

which contradicts the definition of $Q_q$.

(ii) By Proposition 2.4, $(\widehat{\beta}_q, \widehat{Y}_q)$ defines some $\widehat{Q} := Q^{(\widehat{\beta}_q, \widehat{Y}_q)} \in \mathcal{Q}^q \cap \overline{\mathcal{Q}}$ and by Corollary 2.3, we have $f^q(\widehat{Q}|P) = e^{Tk_q(\widehat{\beta}_q, \widehat{Y}_q)}$. Suppose that there exists $Q' = Q^{(\beta', Y')} \in \mathcal{Q}^q$ such that

(2.4) $$f^q(Q'|P) < f^q(\widehat{Q}|P).$$

By 1 of Theorem 2.6, we may assume that $Q' \in \overline{\mathcal{Q}}$, so by Proposition A.7, we have $(\beta', Y') \in \mathcal{A}_q$ and by Corollary 2.3, $f^q(Q'|P) = e^{Tk_q(\beta', Y')}$. However, (2.4) then implies that $k_q(\beta', Y') < k_q(\widehat{\beta}_q, \widehat{Y}_q)$, a contradiction to $(\widehat{\beta}_q, \widehat{Y}_q)$ solving $(\mathcal{P}_q)$. Consequently, $\widehat{Q} = Q_q$. □

REMARK 2.8. Existence results for $Q_q$ can also be found in Bellini and Frittelli [2], but they do not give an explicit formula for $Q_q$; also see the related article [15] of Kramkov and Schachermayer.

2.2. *Sufficient conditions.* In this section, we introduce and discuss conditions $(\mathcal{C}_q)$ which are sufficient for the existence of a solution to $(\mathcal{P}_q)$. In particular, $(\mathcal{C}_q)$ yields an explicit expression for this solution and hence also for the Girsanov parameters of $Q_q$. Although $(\mathcal{C}_q)$ is, in general, only a sufficient condition, it has the advantage that in explicit models it is usually easier to check $(\mathcal{C}_q)$ than to find a solution to $(\mathcal{P}_q)$. We now deduce $(\mathcal{C}_q)$ via a formal application of the Kuhn–Tucker theorem to $(\mathcal{P}_q)$; see Theorem 28.3 in Rockafellar [19]. It then remains to check that $(\mathcal{C}_q)$ implies a solution to $(\mathcal{P}_q)$.

($\mathcal{C}_q$): There exists $\widetilde{\lambda}_q \in \mathbb{R}^d$ with $\widetilde{Y}_q(x) := ((q-1)\widetilde{\lambda}_q^* x + 1)^{1/(q-1)} > 0$, $K$-a.e. and such that ($\mathcal{M}$) is satisfied for $(\widetilde{\beta}_q, \widetilde{Y}_q)$, where $\widetilde{\beta}_q := \widetilde{\lambda}_q$.

THEOREM 2.9. Let $\mathcal{Q}^q \neq \varnothing$ and let ($\mathcal{C}_q$) hold. $(\widetilde{\beta}_q, \widetilde{Y}_q)$ then solves ($\mathcal{P}_q$), that is, $(\widetilde{\beta}_q, \widetilde{Y}_q) = (\widehat{\beta}_q, \widehat{Y}_q)$.

PROOF. We show that $(\widetilde{\beta}_q, \widetilde{Y}_q)$ solves ($\mathcal{P}_q$). First, note that $\mathcal{A}_q \neq \varnothing$ since $\mathcal{Q}^q \neq \varnothing$. In fact, the latter implies, by Theorem 2.6, that $\mathcal{Q}^q \cap \overline{\mathcal{Q}} \neq \varnothing$ and $\mathcal{A}_q \neq \varnothing$ follows from Corollary 2.3 and Proposition A.7 in Appendix A. Let $(\beta, Y) \in \mathcal{A}_q$. Convexity and the definition of $(\widetilde{\beta}_q, \widetilde{Y}_q)$ imply that

$$(2.5) \qquad \beta^* c \beta - \widetilde{\beta}_q^* c \widetilde{\beta}_q \geq 2\widetilde{\beta}_q^* c (\beta - \widetilde{\beta}_q) = 2\widetilde{\lambda}_q^* c (\beta - \widetilde{\beta}_q)$$

and

$$\begin{aligned} g_q(Y(x)) - g_q(\widetilde{Y}_q(x)) &\geq g_q'(\widetilde{Y}_q(x))(Y(x) - \widetilde{Y}_q(x)) \\ &= q(\widetilde{Y}_q^{q-1} - 1)(Y(x) - \widetilde{Y}_q(x)) \\ &= q(q-1)\widetilde{\lambda}_q^* x (Y(x) - \widetilde{Y}_q(x)). \end{aligned}$$

Next, we average both sides of the previous expression with respect to $K$; note that everything is well defined (possibly $-\infty$) since $k_q(\beta, Y) < \infty$ and $g_q(\cdot) \geq 0$ on $(0, \infty)$. Since $(\beta, Y)$ and $(\widetilde{\beta}_q, \widetilde{Y}_q)$ satisfy ($\mathcal{M}$), we thus obtain

$$\int_{\mathbb{R}^d} (g_q(Y(x)) - g_q(\widetilde{Y}_q(x))) K(dx) \geq q(q-1) \widetilde{\lambda}_q^* c (\widetilde{\beta}_q - \beta).$$

This, together with (2.5) and $k_q(\beta, Y) < \infty$, implies that $k_q(\beta, Y) \geq k_q(\widetilde{\beta}_q, \widetilde{Y}_q)$ and the proof is complete. $\square$

REMARK 2.10. (i) In ($\mathcal{C}_q$), we assumed $\widetilde{Y}_q(x) > 0$, $K$-a.e. but not that $\widetilde{Y}_q(x) \geq 0$ for all $x \in \mathbb{R}^d$, which is presumed for $(\widetilde{\beta}_q, \widetilde{Y}_q)$ to be a solution to ($\mathcal{P}_q$). However, we identify all functionals on $\mathbb{R}^d$ which are $K$-a.e. equal since they will describe the same probability measure.

(ii) The assumption $\widetilde{Y}_q(x) > 0$, $K$-a.e. in ($\mathcal{C}_q$) looks very restrictive since it holds if and only if $(q-1)\widetilde{\lambda}_q^* x > -1$, $K$-a.e. However, we assumed that $S = \mathcal{E}(L)$ is strictly positive, that is, that $\Delta L^i > -1$ for $i \in \{1, \ldots, d\}$, so that $\text{supp}(K) \subseteq [-1, \infty)^d$.

(iii) In Theorem 2.9, $\mathcal{Q}^q \neq \varnothing$ can be replaced by

$$(2.6) \qquad \int_{\mathbb{R}^d} g_q(\widetilde{Y}_q(x)) K(dx) < \infty.$$

In fact, $\mathcal{Q}^q \neq \varnothing$ enters only via $\mathcal{A}_q \neq \varnothing$. However, (2.6) implies that $(\widetilde{\beta}_q, \widetilde{Y}_q) \in \mathcal{A}_q$, that is, that $\mathcal{A}_q \neq \varnothing$.

(iv) For $q = 2$, condition (2.6) holds if and only if $\int (\widetilde{\lambda}_2^* x)^2 K(dx) < \infty$. In particular, if $d = 1$ (and $\widetilde{\lambda}_2 \neq 0$), this is equivalent to (local) $P$-square integrability of $L$; see Proposition II.2.29 of Jacod and Shiryaev [13].

Condition $(\mathcal{C}_q)$ appears for $q \in (-\infty, 0)$ also in Kallsen [14] and in Goll and Rüschendorf [10]. However, they do not explain their motivation for the definition of $(\mathcal{C}_q)$. Instead, they prove by direct calculation that $(\mathcal{C}_q)$ determines the qMMM since the conditions which yield an optimal strategy for power utility (resp. an $f$-projection) are satisfied. In a very recent work, Choulli, Stricker and Li [4] state $(\mathcal{C}_q)$ for general $q$, but for the minimal Hellinger martingale measure $\mathcal{Q}$ of order $q$ instead of $Q_q$. We finish this section by linking $\mathcal{Q}$ and $Q_q$ in the present Lévy setting. This is not only interesting per se, but gives a better understanding of $(\mathcal{C}_q)$. For the entire discussion, we assume that $\mathcal{Q}^q \neq \varnothing$. By Propositions 3.3 and 3.5 in Choulli, Stricker and Li [4], $\mathcal{Q} = Q^{(q\beta, qY)}$ can be characterized as that local martingale measure for $S$ (or $L$) such that $\int k_q(\beta_t, Y_t)\, dt - \int k_q(^q\beta_t, {}^q Y_t)\, dt$ is an increasing process for all local martingale measures $Q = Q^{(\beta,Y)}$ with $dQ/dP = \mathcal{E}(N)$ and $N$ locally $q$-integrable. Applying analogous arguments as for Theorem 2.6, one can show that it suffices to consider $\mathcal{Q}^q \cap \overline{\mathcal{Q}}$ and conclude that (the Girsanov parameters of) $\mathcal{Q}_q$ can also be obtained as the solution to $(\mathcal{P}_q)$. Consequently, $\mathcal{Q}$ and $\mathcal{Q}_q$ correspond in the present Lévy setting, provided they exist. Thus, it is no surprise that the sufficient condition given in [4] is also very similar to $(\mathcal{C}_q)$. More precisely, with $\mathcal{D} := \{\lambda \in \mathbb{R}^d | 1 + \inf_{x \in \mathrm{supp}(K)} \lambda^* x \geq 0\}$, their condition reads $\mathrm{int}(\mathcal{D}) \neq \varnothing$ and

$$(2.7) \qquad \int_{\mathbb{R}^d} |(1 + \lambda^* x)^{1/(q-1)} x - h(x)| K(dx) < \infty \qquad \text{for all } \lambda \in \mathrm{int}(\mathcal{D}).$$

However, they use (2.7) only to ensure the existence of some $\lambda \in \mathrm{int}(\mathcal{D})$ satisfying $(\mathcal{M})$; see Lemma 4.4 in [4]. Moreover, in Corollary 4.6, they show that if $q > 1$ and $\mathrm{int}(\mathcal{D}) \neq \varnothing$, then $S$ (or equivalently $L$) is locally $p$-integrable for $p := q/(q-1)$ if and only if (2.7) holds. Consequently, if $q > 1$ and $S$ is locally $p$-integrable, then $(\mathcal{C}_q)$ holds as soon as $\mathrm{int}(\mathcal{D}) \neq \varnothing$.

**3. The variance minimal martingale measure.** In this section, we study the important special case $q = 2$ and refer to $Q_2$ as the *variance minimal martingale measure* (VMMM). We relate $Q_2$ and $(\mathcal{C}_2)$ to the *variance optimal signed martingale measure* (VOSMM) $\widetilde{P}$ and the so-called *structure condition* (SC). Moreover, we discuss in which cases the (easily checked) condition $(\mathcal{C}_2)$ is not only sufficient but also necessary for the existence of $Q_2$.

3.1. *Connection to the variance optimal signed martingale measure.* The VMMM $Q_2$ is obviously related to the *variance optimal signed martingale measure* $\widetilde{P}$ which arises in the mean-variance hedging approach; see Schweizer [21] for an overview and terminology not explained here. The measure $\widetilde{P}$ is obtained by minimizing the variance of the density $dQ/dP$ over all *signed* local martingale measures $Q$ for $S$. Therefore, $Q_2$ and $\widetilde{P}$ coincide if and only if $\widetilde{P}$ is a probability measure equivalent to $P$. In a Lévy setting, $\widetilde{P}$ corresponds to the *minimal signed martingale measure* $\widehat{P}$ which appears in the local risk minimizing hedging approach. It is well known that $\widehat{P}$ exists (even in a more general setting) if the structure condition (SC) is satisfied; (SC) also provides an explicit formula for the density of $\widehat{P}$. Consequently, there should be a link between (SC) and $(\mathcal{C}_2)$. The latter is discussed in this section.

We first introduce (SC) in the Lévy setting. For this, we need $L$ to be a special semimartingale so that by Jacod and Shiryaev ([13], Corollary II.2.38 and Proposition II.2.29), we have

$$L_t = (L_t^c + x * (\mu^L - \nu^P)_t) + \left(b + \int_{\mathbb{R}^d}(x - h(x))K(dx)\right)t$$
$$=: M_t + \gamma t =: M_t + A_t.$$

In addition, (SC) requires that the (local) martingale $M$ is locally $P$-square integrable, that is, $\int \|x\|^2 K(dx) < \infty$. $M$ is then $P$-square integrable and

$$\langle M \rangle_t = \left(c + \int_{\mathbb{R}^d} xx^* K(dx)\right)t =: \sigma t.$$

(SC) is satisfied if there exists a $d$-dimensional predictable process $\lambda$ with

$$(3.1) \qquad A = \int d\langle M\rangle \lambda \quad \text{and} \quad \widehat{K}_T := \int \lambda^* d\langle M\rangle \lambda < \infty;$$

see Definition 1.1 in Choulli and Stricker [3] and Section 12.3 in Delbaen and Schachermayer [5] for a related discussion. Since $A_t = \gamma t$ and $\langle M\rangle_t = \sigma t$, (3.1) holds if and only if there exists $\lambda \in \mathbb{R}^d$ such that $\gamma = \sigma\lambda$ or, equivalently, such that

$$(3.2) \qquad b + \int_{\mathbb{R}^d}(x - h(x))K(dx) = \left(c + \int_{\mathbb{R}^d} xx^* K(dx)\right)\lambda.$$

Under (SC), we can define $\widehat{N} := -\int \lambda^* dM$. If $\widehat{Z} := \mathcal{E}(\widehat{N})$ is a $P$-martingale, then $\frac{d\widehat{P}}{dP} := \widehat{Z}_T$ defines a signed measure called the *minimal signed martingale measure* for $L$. By Proposition 2 of Schweizer [20], it is a local martingale measure for $L$ in the sense that $\widehat{Z}L$ is a local $P$-martingale. Note that if $\widehat{Z} > 0$, that is, if $-\lambda^*\Delta M > -1$, then $\widehat{Z} = \mathcal{E}(\widehat{N})$ is a local martingale and, as in the proof of Proposition A.6, an application of Theorem II.5 of Lépingle and Mémin [17] yields that it is a $P$-martingale, that is, $\widehat{P} \in \mathcal{M}^e(L)$. Since the *mean-variance trade-off process*

$$\widehat{K}_t := \int_0^t \lambda^* \, dA_s = \left\langle \int \lambda^* \, dM \right\rangle_t = \lambda^* \sigma \lambda t$$

is deterministic, Theorem 8 of Schweizer [20] implies that $\widehat{P}$ is equal to the variance-optimal signed martingale measure $\widetilde{P}$. The measures $\widetilde{P}$ and $Q_2$ coincide if $\widetilde{P} = \widehat{P}$ is an equivalent probability measure, that is, if $\widehat{Z}_T > 0$.

We next show that if $\int \|x\|^2 K(dx) < \infty$, then (SC) together with $\widehat{Z}_T > 0$ is equivalent to $(\mathcal{C}_2)$. In fact, (SC) with $\widehat{Z}_T > 0$ holds by (3.2) if and only if there exists $\lambda \in \mathbb{R}^d$ such that

$$Y(x) := -\lambda^* x + 1 > 0, \qquad K\text{-a.e.},$$

$$\int_{\mathbb{R}^d} |x(-\lambda^* x + 1) - h(x)| K(dx) < \infty,$$

$$b - c\lambda + \int_{\mathbb{R}^d} (x(-\lambda^* x + 1) - h(x)) K(dx) = 0.$$

With the replacement $\widetilde{\lambda}_2 := -\lambda$, this equals $(\mathcal{C}_2)$; note that the assumption $\int \|x\|^2 K(dx) < \infty$ implies (2.6). By Theorem 1 of Schweizer [20], if either $L$ is continuous and $\mathcal{M}^e(S) \neq \varnothing$, or $\int \|x\|^2 K(dx) < \infty$ and $\mathcal{Q}^2 \neq \varnothing$, then (SC) holds and $\widehat{P}$ exists. In the first case, we have $\widehat{Z}_T > 0$, so $(\mathcal{C}_2)$ equals (SC); in that case, $Q_2$ always exists. In the second case, $(\mathcal{C}_2)$ reduces to $\widehat{Z}_T > 0$; see also the discussion at the end of Section 2.2.

3.2. *Necessary conditions.* In this section, we give conditions under which $(\mathcal{C}_2)$ is also necessary for the existence of $Q_2$. As suggested by part 4 of Remark 2.10, we assume for the rest of the section that $L$ is (locally) $P$-square integrable, that is, we make the following assumption.

ASSUMPTION. $\int_{\mathbb{R}^d} \|x\|^2 K(dx) < \infty$. Recall that $(b, c, K)$ denotes the characteristic triplet of $L$. To state the main theorem of this section, we introduce the following condition.

$(\mathcal{D})$: *The Girsanov parameters* $(\beta_2, Y_2)$ *of* $Q_2 = Q^{(\beta_2, Y_2)}$ *are such that*

$$\mathcal{H}^1(Y_2) := \left\{ \Psi \in L^2(K) \Big| \int \Psi(x) x K(dx) = 0 \in \mathbb{R}^d, \right.$$

$$\left. |\Psi(x)| \leq a Y_2(x) \ K\text{-a.e. for some } a > 0 \right\}$$

*is* $L^2(K)$-*dense in*

$$\mathcal{H}^2 := \left\{ \Psi \in L^2(K) \Big| \int \Psi(x) x K(dx) = 0 \in \mathbb{R}^d \right\};$$

*elements of $L^2(K)$ are $\mathbb{R}$-valued.*

THEOREM 3.1. *If the VMMM $Q_2 = Q^{(\beta_2, Y_2)}$ exists in $\mathcal{Q}^2$, then $(\mathcal{C}_2)$ is satisfied for some $\widetilde{\lambda}_2$ and $(\widetilde{\lambda}_2, \widetilde{Y}_2) = (\beta_2, Y_2)$ in both of the following cases:*

1. *$c$ is invertible;*
2. *$c = 0$ and $(\mathcal{D})$ holds.*

REMARK 3.2. By Theorem 2.7, we could equivalently state condition $(\mathcal{D})$ and Theorem 3.1 in terms of the solution $(\widehat{\beta}_2, \widehat{Y}_2)$ to $(\mathcal{P}_2)$ instead of the Girsanov parameters $(\beta_2, Y_2)$ of $Q_2 = Q^{(\beta_2, Y_2)}$.

Before we prove Theorem 3.1, we show that if $c = 0$ and $d = 1$, that is, if $L$ is one-dimensional and has no Brownian part, then condition $(\mathcal{D})$ is automatically satisfied.

LEMMA 3.3. *If $Q_2 = Q^{(\beta_2, Y_2)}$ exists in $\mathcal{Q}^2$, $d = 1$, $c = 0$ and $\mathrm{supp}(K) \neq \varnothing$, then $(\mathcal{D})$ holds true.*

PROOF. For an arbitrary $\Psi \in \mathcal{H}^2$, we construct a sequence $(\Psi_n)_{n \in \mathbb{N}} \subseteq \mathcal{H}^1(Y_2)$ converging to $\Psi$ in $L^2(K)$. To this end, we define, for $n \in \mathbb{N}$,

$$A_n := \{x \in \mathbb{R} | |\Psi(x)| \leq n Y_2(x)\} \quad \text{and} \quad \alpha_n := \int_{\mathbb{R}} \Psi(x) \mathbf{I}_{A_n}(x) x K(dx).$$

By the dominated convergence theorem, $\lim_{n \to \infty} \alpha_n = \int \Psi(x) x K(dx) = 0$; this uses the fact that $\int x^2 K(dx) < \infty$. Set $\delta(x) := \mathrm{sign}(x)(|x| \wedge Y_2(x))$ and note that $|\delta(x)| \leq Y_2(x)$ and that $\delta \in L^2(K)$ implies $\delta(x)|x| \in L^1(K)$. Therefore,

$$\gamma := \int_{\mathbb{R}} \delta(x) x K(dx) = \int_{\mathbb{R}} |x|(|x| \wedge Y_2(x)) K(dx) < \infty$$

and $\gamma > 0$ since $\mathrm{supp}(K) \neq \varnothing$ and $K(\{0\}) = 0$ implies that $x\delta(x) > 0$, $K$-a.e. Let $\Psi_n(x) := \Psi(x) \mathbf{I}_{A_n} - \frac{\alpha_n}{\gamma} \delta(x)$ so that $\Psi_n \in L^2(K)$,

$$|\Psi_n(x)| \leq |\Psi(x) \mathbf{I}_{A_n}| + \frac{|\alpha_n|}{\gamma} |\delta(x)| \leq \left(n + \frac{|\alpha_n|}{\gamma}\right) Y_2(x)$$

and

$$\int_{\mathbb{R}} \Psi_n(x) x K(dx) = \alpha_n - \frac{\alpha_n}{\gamma} \int_{\mathbb{R}} \delta(x) x K(dx) = 0,$$

so that $\Psi_n \in \mathcal{H}^1(Y_2)$. Moreover, by the dominated convergence theorem and since $\lim_{n \to \infty} \alpha_n = 0$, we have

$$\lim_{n \to \infty} \int_{\mathbb{R}} |\Psi(x) - \Psi_n(x)|^2 K(dx)$$

$$\leq \lim_{n \to \infty} 2 \left( \int_{\mathbb{R}} |\Psi(x) \mathbf{I}_{A_n^c}|^2 K(dx) + \int_{\mathbb{R}} \left| \frac{\alpha_n}{\gamma} \delta(x) \right|^2 K(dx) \right) = 0.$$

This completes the proof. □

PROOF OF THEOREM 3.1. The main idea is to exploit Theorem 2.7, that is, the fact that $(\beta_2, Y_2)$ solves $(\mathcal{P}_2)$, and deduce that under the assumptions of Theorem 3.1, the solution of $(\mathcal{P}_2)$ has the form as in $(\mathcal{C}_2)$. First, we construct elements in $(\mathcal{A}_2)$ close to $(\beta_2, Y_2)$. Fix $\phi \in \mathbb{R}^d$ and $\Psi \in L^2(K)$ with $|\Psi(x)| \leq aY_2(x)$, $K$-a.e. for some $a > 0$ and

$$(3.3) \qquad c\phi + \int_{\mathbb{R}^d} \Psi(x) x K(dx) = 0;$$

if $c = 0$, this means that $\Psi \in \mathcal{H}^1(Y_2)$. For $\varepsilon > 0$ sufficient small, $(\beta_2 + \varepsilon\phi, Y_2 + \varepsilon\Psi)$ is in $\mathcal{A}_2$; $\Psi$ can be modified on a set of $K$-measure zero. Since $g_2(y) = (y-1)^2$ and $\int g_2(Y_2(x))K(dx) < \infty$, we can define

$$H_{\phi,\Psi}(\varepsilon) := (\beta_2 + \varepsilon\phi)^* c(\beta_2 + \varepsilon\phi) + \int_{\mathbb{R}^d} (Y_2(x) + \varepsilon\Psi(x) - 1)^2 K(dx)$$

and obtain

$$\frac{d}{d\varepsilon} H_{\phi,\Psi}(\varepsilon) = 2\bigg(\varepsilon\phi^* c\phi + \beta_2^* c\phi + \int_{\mathbb{R}^d} \Psi(x)(Y_2(x) + \varepsilon\Psi(x) - 1) K(dx)\bigg).$$

Since $(\beta_2, Y_2)$ solves $(\mathcal{P}_2)$, we have

$$(3.4) \qquad 0 = \frac{d}{d\varepsilon} H_{\phi,\Psi}(0) = 2\bigg(\beta_2^* c\phi + \int_{\mathbb{R}^d} \Psi(x)(Y_2(x) - 1) K(dx)\bigg).$$

We now proceed separately for the two cases of Theorem 3.1, starting with 1. Here, (3.4), together with (3.3), yields

$$(3.5) \qquad \int_{\mathbb{R}^d} (\beta_2^* x - (Y_2(x) - 1))\Psi(x) K(dx) = 0.$$

Moreover, $\Psi \in L^2(K)$ can be chosen arbitrarily under the condition that $|\Psi(x)| \leq aY_2(x)$, $K$-a.e., for some $a > 0$ because (3.3) can always be satisfied by setting $\phi := -c^{-1} \int x\Psi(x) K(dx)$. Consequently, (3.5) implies that

$$(3.6) \qquad \beta_2^* x - (Y_2(x) - 1) = 0, \qquad K\text{-a.e.}$$

In fact, suppose that $\beta_2^* x - (Y_2(x) - 1) > 0$ on a set $A \subseteq \mathbb{R}^d$ with $K(A) > 0$. Then $\widetilde{\Psi}(x) := (\sqrt{\|x\|^2 \wedge 1} \wedge Y_2(x)) \mathbf{I}_A \in L^2(K)$, $\widetilde{\Psi} > 0$, $K$-a.e., $|\widetilde{\Psi}| \leq Y_2$ and

$$\int_{\mathbb{R}^d} (\beta_2^* x - (Y_2(x) - 1))\widetilde{\Psi}(x) K(dx) > 0$$

contradicts (3.5). Since $(\beta_2, Y_2)$ solves $(\mathcal{P}_2)$, we thus obtain from (3.6) that $\widetilde{\lambda}_2 := \beta_2$ satisfies $(\mathcal{C}_2)$. This proves part 1. For part 2, we introduce $\mathcal{L}^0 \subseteq L^2(K)$ as the subspace of all (equivalence classes of) linear functionals. Since $\mathcal{L}^0$ is $d$-dimensional, it is closed and $\mathcal{L}^0 = (\mathcal{L}^{0\perp})^\perp = \mathcal{H}^{2\perp}$. Therefore, (3.4) with $c = 0$ and $(\mathcal{D})$ yield $(Y_2(x) - 1) \in \mathcal{H}^{2\perp} = \mathcal{L}^0$ and hence $(Y_2(x) - 1) = \alpha^* x$, $K$-a.e. for $\alpha \in \mathbb{R}^d$. Setting $\widetilde{\lambda}_2 := \alpha$ implies 2 since $(\beta_2, Y_2)$ solves $(\mathcal{P}_2)$. □

**4. Convergence to the minimal entropy martingale measure.** In this section, we discuss the relationship between $Q_q$ and the *minimal entropy martingale measure* (MEMM) $P_e$ which minimizes the relative entropy $H(Q|P) = E[dQ/dP \log dQ/dP]$ over all local martingale measures $Q$ for $S$. More precisely, we show that under some technical assumptions, $Q_q$ converges for $q \searrow 1$ in entropy to $P_e$; this extends a result of Grandits and Rheinländer from continuous to Lévy processes. In particular, we prove convergence of the Girsanov parameters. At the end of this section, we give a general example in which all assumptions we impose are satisfied.

We first describe $Q_q$ and $P_e$ via roots of functions. Assume that $(\mathcal{C}_q)$ holds for all $q \in (1, 1+\varepsilon]$ for some $\varepsilon > 0$ and that either $\mathcal{Q}^{q'} \neq \varnothing$ or (2.6) is satisfied for $q' := 1+\varepsilon$. Then, for $q \in (1, 1+\varepsilon]$, there exists $\widetilde{\lambda}_q \in \mathbb{R}^d$ such that

$$(4.1) \qquad (q-1)\widetilde{\lambda}_q^* x + 1 > 0, \qquad K\text{-a.e.},$$

$\int |x((q-1)\widetilde{\lambda}_q^* x + 1)^{1/(q-1)} - h(x)| K(dx) < \infty$ and $\Phi(\widetilde{\lambda}_q, q) = 0$ with

$$(4.2) \qquad \Phi(\lambda, q) := b + c\lambda + \int_{\mathbb{R}^d} (x((q-1)\lambda^* x + 1)^{1/(q-1)} - h(x)) K(dx).$$

By Theorem 2.9 and Remark 2.10, $Q_q = Q^{(\beta_q, Y_q)}$ exists for $q \in (1, q']$ with $\beta_q = \widetilde{\lambda}_q$ and $Y_q(x) := ((q-1)\widetilde{\lambda}_q^* x + 1)^{1/(q-1)}$. A similar existence criterion is known for $P_e$; see Theorem 3.1 of Fujiwara and Miyahara [8] or Esche and Schweizer [6], Theorem B and Lemma 15. In fact, if $\lambda_e \in \mathbb{R}^d$ exists with $\int |xe^{\lambda_e^* x} - h(x)| K(dx) < \infty$ and $\Phi_e(\lambda_e) = 0$ for

$$\Phi_e(\lambda) := b + c\lambda + \int_{\mathbb{R}^d} (xe^{\lambda^* x} - h(x)) K(dx),$$

then $P_e$ exists and its Girsanov parameters are

$$\beta_e := \lambda_e \quad \text{and} \quad Y_e(x) := e^{\lambda_e^* x}.$$

For any $\lambda \in \mathbb{R}^d$, we have $\lim_{q \searrow 1}((q-1)\lambda^* x + 1)^{1/(q-1)} = e^{\lambda^* x}$. If, in addition, $(q'-1)\lambda^* x + 1 > 0$, $K$-a.e., and sufficient integrability conditions hold, then $\lim_{q \searrow 1} \Phi(\lambda, q) = \Phi_e(\lambda)$. Consequently, it is natural to expect that the solutions $\widetilde{\lambda}_q$ to $\Phi(\lambda, q) = 0$ also converge to the solution $\lambda_e$ to $\Phi_e(\lambda) = 0$. We show this by an application of the implicit function theorem. Therefore we further assume the existence of some open set $G \subseteq \mathbb{R}^{d+1}$ such that $\{(q, \widetilde{\lambda}_q) | q \in (1, 1+\varepsilon)\} \subseteq G$, $\Phi$ is well defined and continuously differentiable on $G$ and $\det(\frac{d}{d\lambda}\Phi(q, \widetilde{\lambda}_q)) \neq 0$ for all $q \in (1, 1+\varepsilon)$. There then exists a continuously differentiable function $\lambda(q)$ defined on $(1, 1+\varepsilon)$ such that $\lambda(q) = \widetilde{\lambda}_q$ there. Moreover, we assume that $\widetilde{\lambda}_1 := \lim_{q \searrow 1} \lambda(q)$ exists and that $\lim_{q \searrow 1} \Phi(q, \lambda(q)) = \Phi_e(\widetilde{\lambda}_1)$; this holds if $\lambda(\cdot)$ is bounded on $(1, 1+\varepsilon)$ and if we can interchange limit and integration in (4.2). Since $\Phi(q, \lambda(q)) \equiv 0$ we then have $\Phi_e(\widetilde{\lambda}_1) = 0$ and hence $\widetilde{\lambda}_1 = \lambda_e$, as required. Obviously, we also have

$$(4.3) \qquad \lim_{q \searrow 1} \beta_q = \lim_{q \searrow 1} \lambda(q) = \lambda_e = \beta_e,$$

$$(4.4) \qquad \lim_{q \searrow 1} Y_q(x) = \lim_{q \searrow 1}((q-1)\lambda^*(q)x + 1)^{1/(q-1)} = e^{\lambda_e^* x} = Y_e(x), \qquad K\text{-a.e.}$$

In particular, the Lévy measure of $Q_q$ converges to that of $P_e$, that is,

$$\lim_{q \searrow 1} K^{\text{qMMM}}(dx) = \lim_{q \searrow 1} Y_q(x) K(dx) = Y_e(x) K(dx) = K^{\text{MEMM}}(dx);$$

see Proposition A.1 below. Finally, we show that $Q_q$ converges to $P_e$ in entropy, that is, that the relative entropy of $Q_q$ with respect to $P_e$,

$$H(Q_q|P_e) = E_{Q_q}\left[\log \frac{dQ_q}{dP_e}\right] = E_{Q_q}[\log Z_T^{Q_q} - \log Z_T^{P_e}],$$

converges to 0 if $q$ decreases to 1. From Proposition A.3, the formula for the stochastic exponential and Proposition II.1.28 of Jacod and Shiryaev [13] together with Lemma B.4, 3 of Theorem 2.6 and Corollary 2.3, we obtain

$$\begin{aligned}
\log Z_T^{Q_q} &= \beta_q^* L_T^c + (Y_q - 1) * (\mu^L - \nu^P)_T \\
&\quad - \tfrac{1}{2} T \beta_q^* c \beta_q + (\log Y_q - (Y_q - 1)) * \mu_T^L \\
&= \beta_q^* L_T^c + (\log Y_q) * (\mu^L - \nu^P)_T \\
&\quad - \tfrac{1}{2} T \beta_q^* c \beta_q + (\log Y_q - (Y_q - 1)) * \nu_T^P.
\end{aligned}$$

Applying the same arguments and replacing Theorem 2.6 and Corollary 2.3 by Esche and Schweizer [6], Theorem A and Lemma 12, we obtain for $\log Z_T^{P_e}$ the same expression with $(\beta_q, Y_q)$ replaced by $(\beta_e, Y_e)$. Thus,

$$\begin{aligned}
\log \frac{dQ_q}{dP_e} &= (\beta_q^* - \beta_e^*) L_T^c + (\log Y_q - \log Y_e) * (\mu^L - \nu^P)_T \\
&\quad - \frac{1}{2} T(\beta_q^* c \beta_q - \beta_e^* c \beta_e) + (\log Y_q - Y_q - \log Y_e + Y_e) * \nu_T^P.
\end{aligned}$$

Recall from Girsanov's theorem that $\nu^{Q_q} := Y_q \nu^P$ is the $Q_q$-compensator of $\mu^L$ and that $\widetilde{L}$ with $\widetilde{L}_t := L_t^c - tc\beta_q$ is a $Q_q$-martingale. If $(\log Y_q - \log Y_e) * \nu^{Q_q}$ is the $Q_q$-compensator of $(\log Y_q - \log Y_e) * \mu^L$, then

$$\begin{aligned}
\log \frac{dQ_q}{dP_e} &= (\beta_q^* - \beta_e^*) \widetilde{L}_T + (\log Y_q - \log Y_e) * (\mu^L - \nu^{Q_q})_T \\
&\quad - \frac{1}{2} T(\beta_q^* c \beta_q - \beta_e^* c \beta_e) + T(\beta_q^* - \beta_e^*) c \beta_q \\
&\quad + ((\log Y_q - \log Y_e) Y_q - (Y_q - Y_e)) * \nu_T^P;
\end{aligned}$$

see Proposition II.1.28 of Jacod and Shiryaev [13]. The first two terms on the right-hand side are local $Q_q$-martingales and $Q_q$-Lévy processes and thus $Q_q$-martingales, so

$$H(Q_q|P_e) = -\tfrac{1}{2}T(\beta_q^* c\beta_q - \beta_e^* c\beta_e) + (\beta_q^* - \beta_e^*)c\beta_q T$$
$$+ E_{Q_q}[((\log Y_q - \log Y_e)Y_q - (Y_q - Y_e)) * \nu_T^P].$$

Thus, if integration and limit are interchangeable, (4.3) and (4.4) imply that

$$\lim_{q \searrow 1} H(Q_q|P_e) = 0.$$

EXAMPLE. We finish with an example which satisfies all assumptions of this section. Assume that $K$ is of the form $K(dx) = f(x)\,dx$, where $f$ is a bounded density such that $\mathrm{supp}(K) \subseteq (-1,\ell)$ with $0 < \ell < \infty$. Moreover, let $L$ be of dimension one and if $L$ has no Brownian part, that is, if $c = 0$, then let it have jumps of positive and negative heights, that is, $K((-1,0)) > 0$ and $K((0,\ell)) > 0$. We show that there exists $\varepsilon > 0$ such that $(\mathcal{C}_q)$ has a solution $\widetilde{\lambda}_q$ for all $q \in (1, 1+\varepsilon]$ and such that we can take

$$G := G(\varepsilon) := \{(q,\lambda)|q \in (1,1+\varepsilon), q \in (\gamma_1(q), \gamma_2(q))\}$$

with $\gamma_1(q) := -\frac{1}{\ell(q-1)}$ and $\gamma_2(q) := \frac{1}{q-1}$. All integrability conditions of Section 4 are then satisfied due to boundedness of $f$ and of $\mathrm{supp}(K)$; note that $y \mapsto y \log y$ is bounded from below and that $Y^e(x) = \exp(\lambda_e^* x)$ is $K$-a.e. bounded away from 0. For $q \in (1,2)$ and $\lambda \in (\gamma_1(q), \gamma_2(q))$, condition (4.1) is satisfied and

$$\frac{d}{d\lambda}\Phi(\lambda, q) = c + \int_{-1}^{\ell}(x^2((q-1)\lambda x + 1)^{(2-q)/(q-1)})f(x)\,dx$$
$$\geq \begin{cases} c + \int_{-1}^{0} x^2 f(x)\,dx =: \delta_1 > 0, & \text{if } \lambda \in (\gamma_1(q), 0], \\ c + \int_{0}^{\ell} x^2 f(x)\,dx =: \delta_2 > 0, & \text{if } \lambda \in [0, \gamma_2(q)). \end{cases}$$

Thus, $\frac{d}{d\lambda}\Phi(\lambda, q) \geq \delta := \min\{\delta_1, \delta_2\} > 0$. Let

$$b_0 := \Phi(\lambda, 0) = \Phi_e(\lambda) = b + \int_{-1}^{\ell}(x - h(x))f(x)\,d(x)$$

and note that $\lim_{q \searrow 1} \gamma_1(q) = -\infty$ and $\lim_{q \searrow 1} \gamma_2(q) = \infty$. If $b_0 < 0$, we can hence find $\varepsilon_1 > 0$ such that $\gamma_2(q) > |b_0|/\delta$ for all $q \in (1, 1+\varepsilon_1]$. Then for all $q \in (1, 1+\varepsilon_1)$, there exists a solution $\widetilde{\lambda}_q \in (0, |b_0|/\delta) \subseteq (0, \gamma_2(q))$ to $\Phi(\lambda, q) = 0$ and we can take $G = G(\varepsilon_1)$; note, in addition, that $q \mapsto \widetilde{\lambda}_q$ is bounded. Analogously, if $b_0 > 0$, we select $\varepsilon_2 > 0$ such that $|\gamma_1(q)| > b_0/\delta$ for all $q \in (1, 1+\varepsilon_2]$, which implies for these $q$ the existence of a solution $\widetilde{\lambda}_q \in$ $(-b_0/\delta, 0) \subseteq (\gamma_1(q), 0)$ to $\Phi(\lambda, q) = 0$ and that we can take $G = G(\varepsilon_2)$. Finally, $b_0 = 0$ is a trivial case, since we then have $\Phi(0, q) = 0$ for all $q > 1$ so that $\beta_q = 0$ and $Y_q = 1$, that is, $Q_q = P = P_e$. This concludes the example.

## APPENDIX A: CHANGE OF MEASURE AND LÉVY PROCESSES

In this appendix we gather the required results on changes of measure and Lévy processes. In particular, we give conditions under which two processes are the Girsanov parameters of an equivalent local martingale measure. For unexplained notation, we refer to Jacod and Shiryaev [13].

We fix a filtered probability space $(\Omega, \mathcal{F}, \mathbb{F}, P)$ with $\mathbb{F} = (\mathcal{F}_t)_{0 \leq t \leq T}$ satisfying the usual conditions under $P$ and $\mathcal{F}_0$ trivial. Moreover, we work throughout with the truncation function $h(x) := x \mathbf{I}_{\{\|x\| \leq 1\}}$. By $\mathbf{P}$, we denote the predictable $\sigma$-field on $\Omega \times [0, T]$ and by $(B, C, \nu)$, the $P$-characteristics of the semimartingale $X$ with respect to $h$. As in Proposition II.2.9 of [13], we can and do always choose a version of the form

$$B = \int b\, dA, \qquad C = \int c\, dA,$$

(A.1)
$$\nu(\omega; dx, dt) = K_{\omega,t}(dx)\, dA_t(\omega),$$

where $A$ is a real-valued, predictable, increasing and locally integrable process, $b$ an $\mathbb{R}^d$-valued predictable process, $c$ a predictable process with values in the set of all symmetric nonnegative definite $d \times d$-matrices and $K_{\omega,t}(dx)$ a transition kernel from $(\Omega \times [0, T], \mathbf{P})$ into $(\mathbb{R}^d, \mathcal{B}^d)$ with $K_{\omega,t}(\{0\}) = 0$ and $\int_{\mathbb{R}^d} (1 \wedge \|x\|^2) K_{\omega,t}(dx) \leq 1$ for all $t \leq T$.

We now turn to the description of absolutely continuous probability measures. The following Girsanov-type result shows that any $Q \ll P$ can be described by two parameters $\beta$ and $Y$.

PROPOSITION A.1 (Theorem III.3.24 of Jacod and Shiryaev [13]).   *Let $X$ be a semimartingale with $P$-characteristics $(B^P, C^P, \nu^P)$ and denote by $c$, $A$ the processes from* (A.1). *For any probability measure $Q \ll P$, there exist a $\mathbf{P} \otimes \mathcal{B}^d$-measurable function $Y \geq 0$ and a predictable $\mathbb{R}^d$-valued process $\beta$ satisfying*

$$\|(Y - 1)h\| * \nu_T^P + \int_0^T \|c_s \beta_s\|\, dA_s + \int_0^T \beta_s^* c_s \beta_s\, dA_s < \infty, \qquad Q\text{-}a.s.$$

*and such that the $Q$-characteristics $(B^Q, c^Q, \nu^Q)$ of $X$ are given by*

$$B_t^Q = B_t^P + \int_0^t c_s \beta_s\, dA_s + ((Y - 1)h) * \nu_t^P,$$

$$C_t^Q = C_t^P,$$

$$\nu^Q(dx, dt) = Y_t(x) \nu^P(dx, dt).$$

*We call $\beta$ and $Y$ the* Girsanov parameters *of $Q$ (with respect to $P$, relative to $X$) and write $Q = Q^{(\beta,Y)}$ to emphasize the dependence.*

REMARK A.2. (i) In Proposition A.1, we have $Y(x) > 0 \, dP \otimes dt$-a.e. for $K$-a.e. $x$ if and only if $Q \approx P$.

(ii) Note that $\beta$ and $Y$ are not unique. In fact, $Y$ is unique only $\nu^P$-a.e., and for fixed $c$ and $A$, we have $A$-a.e. uniqueness only for $c\beta$. In the whole article, we fix a process $L$ and express the Girsanov parameters of any $Q \ll P$ relative to $L$. We then identify all versions of Girsanov parameters $(\beta, Y)$ which describe the same $Q$. In particular, if we say that the Girsanov parameters $(\beta, Y)$ of $Q$ are time-independent, we mean that there exists one version with this property.

**A.1. Lévy processes.** Let $Q \approx P$ and $L = (L_t)_{0 \leq t \leq T}$ be an $\mathbb{F}$-adapted stochastic process with RCLL paths and $L_0 = 0$. Then $L$ is called a $(Q, \mathbb{F})$-*Lévy process* if for all $s \leq t \leq T$, the increment $L_t - L_s$ is independent of $\mathcal{F}_s$ under $Q$ and has a distribution which depends on $t - s$ only. Recall that a Lévy process is a Feller process, so that $\mathbb{F}^{L,Q}$, the $Q$-augmentation of the filtration generated by $L$, automatically satisfies the usual conditions under $Q$. If $Q = P$, we sometimes even omit the mention of $P$, that is, refer to $L$ simply as a Lévy process and write $\mathbb{F}^L$. In particular, if $Q = P$ and $\mathbb{F} = \mathbb{F}^L$ for quantities depending on $P$ and $L$, we often do not write this dependence explicitly; this is done, for example, for Girsanov parameters. We will frequently use the fact that for $Q \approx P$, every $(Q, \mathbb{F})$-Lévy process is an $\mathbb{F}$-semimartingale and a $(Q, \mathbb{F})$-martingale if and only if it is a $(Q, \mathbb{F})$-local martingale; see He, Wan and Yan [12], Theorem 11.46. In addition, Lévy processes have the weak predictable representation property; see [13], Theorem III.4.34. This implies an explicit formula for the density process of any $Q \approx P$.

PROPOSITION A.3 (Proposition 3 of Esche and Schweizer [6]). *Let $L$ be a $P$-Lévy process and $\mathbb{F} = \mathbb{F}^L$. If $Q \approx P$ with Girsanov parameters $(\beta, Y)$, then the density process of $Q$ with respect to $P$ is given by $Z^Q = \mathcal{E}(N^Q)$, with*

$$N^Q_t = \int_0^t \beta_s^* \, dL_s^c + (Y - 1) * (\mu^L - \nu^P)_t.$$

REMARK A.4. We frequently use the fact that for $f:(-1, \infty) \to \mathbb{R}$ sufficiently integrable, we have $\sum_{s \leq t} f(\Delta N_s^Q) = f(Y-1) * \mu_t^L$.

It is well known that a Lévy process can be characterized by the particular structure of its characteristics; see Corollary II.4.19 of [13]. In fact, let $Q \approx P$ and $L$ be a $(Q, \mathbb{F})$-semimartingale. $L$ is then a $(Q, \mathbb{F})$-Lévy process if and only if there exists a version of its $Q$-characteristics such that

$$(A.2) \quad B_t^Q(\omega) = b^Q t, \qquad C_t^Q(\omega) = c^Q t, \qquad \nu^Q(\omega; dx, dt) = K^Q(dx) \, dt,$$

where $b^Q \in \mathbb{R}^d$, $c^Q$ is a symmetric nonnegative definite $d \times d$-matrix and $K^Q$ is a positive measure on $\mathbb{R}^d$. We call $(b^Q, c^Q, K^Q)$ the *characteristic triplet* of $L$ (with respect to $Q$). For a $P$-Lévy process, we drop the mention of $P$ and simply write $(b, c, K)$. Moreover, we always use the notation

$$\nu^P(dx, dt) = K(dx) \, dt.$$

As an immediate consequence of Girsanov's theorem and (A.2), we obtain, for any $(P, \mathbb{F})$-Lévy process $L$, the following well-known characterization of the set of all probability measures $Q \approx P$ under which $L$ is a $(Q, \mathbb{F})$-Lévy process.

COROLLARY A.5. *Let $L$ be an $(P, \mathbb{F})$-Lévy process and $Q = Q^{(\beta, Y)} \approx P$. $L$ is then a $(Q, \mathbb{F})$-Lévy process if and only if $\beta$ and $Y(x)$ are $dP \otimes dt$-a.e. time-independent and deterministic for $K$-a.e. $x \in \mathbb{R}^d$.*

**A.2. Change of measure.** So far, we described for any $Q \approx P$ the corresponding Girsanov parameters. We now want to start with arbitrary predictable processes $\beta$ and $Y$ and give conditions under which they define a probability measure $Q \approx P$ and can be identified as the Girsanov parameters of $Q$. We formulate sufficient integrability conditions in terms of the strictly convex function $g_q : (0, \infty) \to \mathbb{R}_+$, $g_q(y) := y^q - 1 - q(y - 1)$, where $q \in I := (-\infty, 0) \cup (1, \infty)$. This function arises in the computation of $f^q(Q|P) = E[(Z_T^Q)^q]$; see Proposition 2.2.

PROPOSITION A.6. *Let $L$ be a $P$-Lévy process with characteristic triplet $(b, c, K)$, $\mathbb{F} = \mathbb{F}^L$, $q \in I$, $\beta$ a predictable process and $Y \geq 0$ a $\mathbf{P} \otimes \mathcal{B}^d$-measurable function with $Y(x) > 0$, $K$-a.e. If*

$$(A.3) \quad \int_{\mathbb{R}^d} g_q(Y_s(x)) K(dx) \leq \text{const}, \qquad dP \otimes dt\text{-a.e.},$$

*then $Y - 1$ is integrable with respect to $\mu^L - \nu^P$. If, in addition,*

$$(A.4) \quad \beta_s^* c \beta_s \leq \text{const}, \qquad dP \otimes dt\text{-a.e.},$$

*then $Z := \mathcal{E}(N)$ with*

$$(A.5) \quad N_t = \int_0^t \beta_s^* \, dL_s^c + (Y - 1) * (\mu^L - \nu^P)_t$$

*is a strictly positive $P$-martingale. In particular, for $dQ/dP := Z_T$, we have $Q = Q^{(\beta, Y)}$, that is, $\beta$ and $Y$ are the Girsanov parameters of $Q$.*

PROOF. The integrability of $Y - 1$ with respect to $\mu^L - \nu^P$ follows from Lemma B.1 together with Theorem II.1.33 d) in [13]. Thus, by (A.4), $N$ is a local martingale and, in addition, quasi-left-continuous, so by Theorem II.5 in Lépingle and Mémin [17], $\mathcal{E}(N)$ is a martingale if the predictable compensator of $\langle N^c \rangle . + \sum_{s \leq .}((\Delta N_s)^2 \wedge |\Delta N|_s)$ is bounded; note that for Theorem II.5 of [17], it suffices for $N$ to be a local martingale. In addition, $\mathcal{E}(N)$ is strictly positive since $Y > 0$ implies that $\Delta N > -1$, so it only remains to show boundedness of the compensator. For $\langle N^c \rangle = \int \beta_t^* c \beta_t \, dt$, which is already the predictable compensator of itself, the claim is trivial by (A.4). The jump term can be rewritten as

$$\text{(A.6)} \qquad \sum_{s \leq t}((\Delta N_s)^2 \wedge |\Delta N_s|) = ((Y-1)^2 \wedge |Y-1|) * \mu_t^L.$$

Since $N$ is, in particular, a special semimartingale, (A.6) defines, by Propositions II.1.28 and II.2.29 a) of [13], a locally integrable process. Also, by Proposition II.1.28, the latter has $(Y-1)^2 \wedge |Y-1| * \nu^P$ as predictable $P$-compensator. This compensator is then bounded thanks to Lemma B.3 and Assumption (A.3). Finally, $Q = Q^{(\beta,Y)}$ holds by Proposition 7 of Esche and Schweizer [6]. This completes the proof. □

We finish this section with a result which gives conditions for the Girsanov parameters $(\beta, Y)$ of $Q = Q^{(\beta,Y)} \approx P$ to imply that $Q$ is a local martingale measure for a $P$-Lévy processes $L$, that is, for $Q \in \mathcal{M}^e(L)$.

PROPOSITION A.7 (Theorem 3.1 of Kunita [16]). *Let $L$ be a $P$-Lévy process with characteristic triplet $(b, c, K)$, $\mathbb{F} = \mathbb{F}^L$ and $Q = Q^{(\beta,Y)} \approx P$. Then $Q \in \mathcal{M}^e(L)$ if and only if*

$$\int_0^T \int_{\mathbb{R}^d} |xY_t(x) - h(x)| K(dx) < \infty, \qquad P\text{-a.e.},$$

$$b + c\beta_t + \int_{\mathbb{R}^d}(xY_t(x) - h(x))K(dx) = 0, \qquad dP \otimes dt\text{-a.e. on } \Omega \times [0, T];$$

*this condition is called the* martingale condition *for $L$.*

## APPENDIX B: AUXILIARY RESULTS

This section contains some simple auxiliary results.

LEMMA B.1. *Fix $q \in I$. There then exists $c = c(q) > 0$ such that*

$$(1 - \sqrt{y})^2 \leq c g_q(y) \qquad \text{for all } y > 0.$$

PROOF. We need to find $c > 0$ such that $f(y) := cg_q(y) - (1 - \sqrt{y})^2$ is nonnegative on $(0, \infty)$. For $q < 0$, we can take $c = -\frac{1}{q}$ since $f(1) = 0$, $\frac{d}{dy}f(y) < 0$ on $(0,1)$ and $\frac{d}{dy}f(y) > 0$ on $(1, \infty)$. For $q > 1$, take $c > \frac{1}{q-1} > \frac{1}{2q(q-1)}$ and define $\overline{y} := (2cq(q-1))^{1/(1/2-q)}$. Calculating $\frac{d^2}{dy^2}f$ yields strict concavity of $f$ on $(0, \overline{y})$ and strict convexity on $(\overline{y}, \infty)$. Moreover, $\overline{y} < 1$, $\frac{d}{dy}f(1) = 0$ and $f(1) = 0$. Since $f$ is continuous, it is thus nonnegative on $(0, \infty)$ if $f(0) \geq 0$, which holds true by the choice of $c$. □

LEMMA B.2. *Fix $q \in I$ and $\overline{y} > 1$. There then exists a constant $C = C(\overline{y}, q) > 0$ such that for all $c \geq C$,*
$$(y-1)^2 \leq cg_q(y) \qquad \text{for all } y \in (0, \overline{y}].$$

PROOF. Define $f(y) := cg_q(y) - (y-1)^2$. For $q < 2$, let $C := \frac{2\overline{y}^{2-1}}{q(q-1)}$. Then $f$ is convex on $(0, \overline{y})$ with minimum in $y = 1$, where $f(1) = 0$, so $f$ is nonnegative on $(0, \overline{y}]$. For $q \geq 2$, set $C := \frac{2}{q}$ so that $\frac{d}{dy}f(y) \leq 0$ on $(0,1)$ and $\frac{d}{dy}f(y) \geq 0$ on $(1, \overline{y}]$. Since $f(1) = 0$, we have $f(y) \geq 0$ on $(0, \overline{y}]$. □

LEMMA B.3. *For $q \in I$, there exists $C = C(q) > 0$ such that*
$$(y-1)^2 \wedge |y-1| \leq Cg_q(y) \qquad \text{for all } y > 0.$$

PROOF. Lemma B.2 with $\overline{y} = 2$ implies the claim for $0 \leq y \leq 2$. For $y > 2$, note that $(y-1)^2 \geq |y-1| = y-1$ and define $f(y) := Cg_q(y) - (y-1)$, with $C \geq \max\{-\frac{2}{q}, \frac{1}{g_q(2)}\}$ for $q < 0$ and $C \geq \max\{\frac{1}{q(2^{q-1}-1)}, \frac{1}{g_q(2)}\}$ for $q > 1$. The function $f$ is then increasing on $[2, \infty)$ and $f(2) \geq 0$. □

LEMMA B.4. *For $q \in I$ and $y > 0$, we have*
$$\log y - (y-1) \leq g_q(y) \quad \text{and} \quad \log y - (y-1) \leq y \log y - (y-1).$$

PROOF. Both $y \mapsto g_q(y) - (\log y - (y-1))$ and $y \mapsto y \log y - (y-1) - (\log y - (y-1))$ are strictly convex functions on $\mathbb{R}_+$. Their unique minimum is in $y = 1$, where they are equal to 0, so they are nonnegative on $\mathbb{R}_+$. □

## APPENDIX C: OMITTED PROOFS

This section contains the proofs omitted in Section 2.

PROOF OF PROPOSITION 2.2. Itô's formula applied to $Z = Z^Q = \mathcal{E}(N)$ yields
$$Z_t^q = 1 + \int_0^t Z_{s-}^q \left( q\, dN_s + \frac{q(q-1)}{2} d\langle N^c \rangle_s \right)$$
$$+ \sum_{s \leq t} Z_{s-}^q ((\Delta N_s + 1)^q - q(\Delta N_s + 1) - 1 + q).$$

Recall from Proposition A.3 the expression for $N$ and note that $\langle N^c \rangle = \int \beta_t^* c \beta_t \, dt$ and $N$ are locally $P$-integrable and, since $Q \in \mathcal{Q}^q$, so is $Z^q$. Thus, we also have $\sum_{s \leq t} Z_{s-}^q ((\Delta N_s + 1)^q - q(\Delta N_s + 1) - 1 + q)$ and

$$\sum_{s \leq t} ((\Delta N_s + 1)^q - q(\Delta N_s + 1) - 1 + q) = (Y^q - 1 - q(Y-1)) * \mu^L$$
$$= g_q(Y) * \mu^L$$

are locally $P$-integrable. Since $g_q$ is nonnegative, Proposition II.1.28 of [13] then implies that the predictable compensator of $g_q(Y) * \mu^L$ is $g_q(Y) * \nu^P$. Moreover,

$$(Y^q - 1 - q(Y-1)) * (\mu^L - \nu^P) + q(Y-1) * (\mu^L - \nu^P) = (Y^q - 1) * (\mu^L - \nu^P)$$

since both sides are local martingales having the same jumps; see Definition II.1.27 in [13]. From this and the formula for $N$ from Proposition A.3, we obtain the canonical decomposition

$$dZ^q = Z_-^q \bigg( q \, dN^c + d((Y^q - 1) * (\mu^L - \nu^P))$$
$$+ \frac{q(q-1)}{2} d\langle N^c \rangle + d((Y^q - 1 - q(Y-1)) * \nu^P) \bigg)$$

(C.1) $\quad = Z_-^q (d\widehat{M} + d\widehat{A})$
$\quad = d\mathcal{E}(\widehat{M} + \widehat{A})$
$\quad = d(\mathcal{E}(\widehat{M})\mathcal{E}(\widehat{A})),$

where the last equality holds by Yor's formula since $\widehat{A}$ is of finite variation and continuous so that $[\widehat{M}, \widehat{A}] \equiv 0$. Moreover, $Q \in \mathcal{Q}^q$ implies that $Z^q$ is a positive submartingale and thus of class $(D)$ since $0 \leq Z_\tau^q \leq E[Z_T^q | \mathcal{F}_\tau]$ for all stopping times $\tau \leq T$. Since $\widehat{A}$ is nonnegative and continuous, we have $\mathcal{E}(\widehat{A}) = e^{\widehat{A}} \geq 1$, so (C.1) implies that $\mathcal{E}(\widehat{M})$ is a local $P$-martingale of class $(D)$ and thus a martingale; this uses the fact that $\mathcal{E}(\widehat{M})$ is (strictly) positive since $\Delta \widehat{M} > -1$ [because $Y(x) > 0$ $K$-a.s. implies that $Y^q(x) - 1 > -1$ $K$-a.s.]. Moreover, (C.1) then implies the $R^p$-integrability of $\mathcal{E}(\widehat{A})$. This completes the proof. $\square$

PROOF OF PROPOSITION 2.4. Propositions A.6 and A.7 imply that $(\overline{\beta}, \overline{Y})$ are the Girsanov parameters of some $\overline{Q} = Q^{(\overline{\beta}, \overline{Y})} \in \mathcal{M}^e(L) = \mathcal{M}^e(S)$. It remains to show that $\overline{Q} \in \mathcal{Q}^q$. This can be done, as in the proof of Proposition 2.2, by an application of Itô's formula to obtain the canonical decomposition and, in particular, that $f^q(Z_T^{\overline{Q}}) = e^{\widehat{A}_T(\overline{Q})} \mathcal{E}(\widehat{M}(\overline{Q}))_T$. The only difference in the proof is the way one obtains the fact that $g_q(\overline{Y}) * \mu^L$ is locally $P$-integrable; this cannot be done as before since we do not know that $\overline{Q} \in \mathcal{Q}^q$.

However, since $g_q$ is nonnegative, we obtain from (2.3) that $g_q(\overline{Y}) * \nu^P$ is locally $P$-integrable and this is, by Proposition II.1.28 of [13], equivalent to local $P$-integrability of $g_q(\overline{Y}) * \mu^L$. Thus, it only remains to show that $f^q(Q|P) = E[e^{\widehat{A}_T(\overline{Q})}\mathcal{E}(\widehat{M}(\overline{Q}))_T] < \infty$. This holds true since $\Delta \widehat{M}(\overline{Q}) > -1$ implies that $\mathcal{E}(\widehat{M}(\overline{Q}))$ is a $P$-supermartingale and since $\widehat{A}_T(\overline{Q})$ is a constant. $\square$

**Acknowledgments.** Part of this work was completed during a visit of M. Jeanblanc to the Graduate School of Economics, Nagoya City University. She expresses her gratitude for the kind invitation and financial support. S. Klöppel wishes to thank M. Schweizer for many helpful discussions and suggestions and Y. Miyahara wants to thank Prof. H. Kunita for many valuable discussions. Thanks also go to an anonymous referee for his comments which led to an improvement of the presentation.

M. Jeanblanc
Equipe d'analyse et probabilités
Université d'Evry Val d'Essonne
Rue du Père Jarlan
91025 Evry Cedex
France
E-mail: monique.jeanblanc@univ-evry.fr

S. Klöppel
Vienna University of Technology
Financial and Actuarial Mathematics
Wiedner Hauptstrasse 8 /105-1 FAM
1040 Vienna
Austria
E-mail: kloeppel@fam.tuwien.ac.at

Y. Miyahara
Graduate School of Economics
Nagoya City University
Mizuhochou Mizuhoku
Nagoya 8501
Japan
E-mail: y-miya@econ.nagoya-cu.ac.jp